\documentclass[10pt]{IEEEtran}
\IEEEoverridecommandlockouts
\usepackage[utf8]{inputenc}
\usepackage{amsmath}
\usepackage{amsthm}
\usepackage{amssymb}
\usepackage{graphicx}
\usepackage{psfrag}
\usepackage{cite}
\usepackage{url}
\usepackage{enumerate}

 \theoremstyle{plain}

 \theoremstyle{definition}
 \newtheorem{defn}{Definition}
 \newtheorem{rem}{Remark}

 \newtheorem{assm}{Assumption}
\newtheorem{stepwise}{Step}

\newcommand{\bb}[0]{\begin{bmatrix}}
\newcommand{\eb}[0]{\end{bmatrix}}
\newcommand{\be}[0]{\begin{equation}}
\newcommand{\ee}[0]{\end{equation}}
\newcommand{\ben}[0]{\begin{equation*}}
\newcommand{\een}[0]{\end{equation*}}

\title{\vspace{18pt} \bf \Large On Output Feedback Adaptive Control}

\author{Zheng Qu, Daniel Wiese, Anuradha M. Annaswamy and Eugene Lavretsky
\thanks{Z. Qu, D. Wiese and A. M. Annaswamy are with the Department
of Mechanical Engineering, Massaschusetts Institute of Technology, Cambridge,
MA, 02139 e-mail: ({mzqu@mit.edu}).}
\thanks{E.~Lavretsky is with the Boeing Company, Huntington Beach CA, 92648.}}

\begin{document}

\title{Squaring-Up Method In the Presence of Transmission Zeros}
\maketitle
\begin{abstract}
This paper presents a method to square up a generic MIMO system that
already possesses transmission zeros. The proposed method is developed
based on and therefore can be incorporated into the existing method
that has been proven effective on a system without transmission zeros.
It has been shown that for the generic system considering here, the
squaring-up problem can be transformed into a state-feedback pole-placement
problem with uncontrollable modes.
\end{abstract}

\section{Introduction}

Square system plays key role in control theory development because
of some unique properties it may possess, such as invertibility\cite{RefWorks:96}
and strictly positive realness\cite{RefWorks:98}. Some control technique
is first developed based on square systems and later on, extended
into a more general plant using square systems as a leverage. Such
extension usually requires a squaring method. By squaring, we particularly
define it as a way to coincide the number of inputs and outputs and
also, make the squared system minimum phase. Minimum phase system
is of particular interest because its inversion is also stable, which
is the necessary conditions of some advanced control design\cite{RefWorks:55}.

Two different squaring methods have been presented in the previous
literature, i.e. squaring-down and squaring-up. The squaring-down
method is first attempted in 1970s\cite{RefWorks:74,RefWorks:68}
and is revealed to be equivalent to pole-placement using output feedback
in the transformed coordinate. Pole-placement using output feedback
has been shown only achievable under some specific conditions and
therefore can be restrictive. On the other hand, the research on squaring-up
method has been sparse until Misra's work in 1990s\cite{RefWorks:70,RefWorks:102}.
It has been shown the squaring-up method is equivalent to pole-placement
using state feedback in the transformed coordinate and therefore is
much more feasible. Actually, pole-placement using state feedback
has been considered a solved problem as long as the controllability
condition in the transformed coordinate is met. In terms of control
design, however, squaring-up method does not prevail over squaring-down
method because it brings pseudo-inputs or outputs into the system
that cannot be used in the physical world.

Recently, the squaring-up method has gained increasing interest in
the new development of adaptive control theory when a minimum phase
system is assumed\cite{RefWorks:55}, or is required as an intermediate
step\cite{RefWorks:90,RefWorks:36}. Some properties the square system
is having, such as strictly positive realness, can be transmitted
to the original system using proper partition. Thus the pseudo-inputs
and outputs are only used in the gain design but never used in the
actual physical action. The final results of the squaring-up treatment
produces a squaring-down matrix, which can be realized in the real
world. One problem arises, however, that when the given system already
possesses transmission zeros, the existing squaring-up method fails.
Previous literature warns the reader but did not specify a solution\cite{RefWorks:70}.
Although the existence of transmission zeros in a non-square system
is especially rare\cite{RefWorks:92}, the failure of squaring-up
method in such case draws interest because the existing transmission
zero, as long as it is nonminimum-phase, has been demonstrated to
be nonpreventive in the adaptive control design. Such observation
motivates us to rationalize the failure and its countermeasure. This
paper will show the pre-existing nonminimum phase zeros are the only
case of our interest that the existing squaring-up method cannot work.
Section \ref{sub:Countermeasure} will provide a remedy to the method
and Section \ref{sec:Example} will present a numeric example.

\section{Preliminaries}

Given a system $\{A,B,C,D\}$, $A\in\mathbb{R}^{n\times n}$, $B\in\mathbb{R}^{n\times m}$,
$C\in\mathbb{R}^{p\times n}$, $D\in\mathbb{R}^{p\times m}$. The
system $\{A,B,C,D\}$ is square if $m=p$. It is tall if $m<p$. It
is fat if $m>p$. For a non-square system, the procedure to square-up
is defined as adding more inputs or outputs to make $m=p$. The procedure
to square-down is defined as abandoning inputs or outputs to make
$m=p$. Given the slight disagreement in the definition of the zeros
of a MIMO system\cite{RefWorks:104,RefWorks:106}, this paper begins
with a series of zero definitions that are widely accepted in recent
literature and relevant to our investigation. 
\begin{defn}
The normal rank of a matrix function $X(s)$ is defined as the rank
of $X(s)$ for almost all the values of $s$.
\end{defn}
It is the interest of this paper to study the few $s$ such that $X(s)$
loses its normal rank.
\begin{defn}
\cite{RefWorks:91}The Rosenbrock matrix of a system $\{A,B,C,D\}$
is defined as: 
\begin{equation}
R(s)=\left[\begin{array}{cc}
sI-A & -B\\
C & D
\end{array}\right]\label{eq:Rosenbrock}
\end{equation}

\end{defn}
The Rosenbrock matrix is first proposed in the reference\cite{RefWorks:91}
and has been widely used to study the zeros of MIMO systems. 
\begin{defn}
If for a system $\{A,B,C,D\}$, the rank of $R(s)$ is strictly less
than the $min(n+m,\ n+p)$ for any complex values of $s$, then the
system is degenerate. 
\end{defn}
It is noted some controllable and observable system can be degenerate.
The system is degenerate when there are some repeated states, inputs
or outputs. For example, if $m<p$ and the system has repeated inputs,
the $B$ has identical columns and naturally, $R(s)$ has a rank less
than $n+m$ for any $s$. 
\begin{defn}
For a system $\{A,B,C,D\}$, the input decoupling zeros are the values
of $s$ such that the following $n\times(n+m)$ matrix loses its normal
rank: 
\begin{equation}
R_{I}(s)=\left[\begin{array}{cc}
sI-A & -B\end{array}\right]
\end{equation}

\end{defn}
The input decoupling zeros are a subset of system poles. They are
actually the uncontrollable modes of the system. 
\begin{defn}
\label{Def:Output Decoup}For a system $\{A,B,C,D\}$, the output
decoupling zeros are the values of $s$ such that the following $(n+p)\times n$
matrix loses its normal rank: 
\begin{equation}
R_{O}(s)=\left[\begin{array}{c}
sI-A\\
C
\end{array}\right]
\end{equation}

\end{defn}
The output decoupling zeros are a subset of system poles. They are
actually the unobservable modes of the system. 
\begin{defn}
For a system $\{A,B,C,D\}$, the invariant zeros are the values of
$s$ such that $R(s)$ loses its normal rank.
\end{defn}
This definition is first proposed by Rosenbrock\cite{RefWorks:91}.
The name ``invariant'' comes from the fact that such zeros are invariant
under feedback action, either state feedback or observer feedback. 
\begin{defn}
\label{Def:Tzero1}For a system $G(s)=C(sI-A)^{-1}B+D$, the transmission
zeros are defined as the values of $s$ such that the rank of $G(s)$
is less than its normal rank.
\end{defn}
The definition is first proposed by Desoer and Schulman and later
on generalized by MacFarlane\cite{RefWorks:107}. It is intended to
describe a special property of the system that it blocks the transmission
from input to output at some specific frequencies. Such frequencies
are named as the ``transmission zeros'' of the system. Recently,
a state-space form of transmission zero definition has been proposed. 
\begin{defn}
\label{Def:Tzero2}\cite{RefWorks:92}For a non-degenerate system
$\{A,B,C,D\}$ that is controllable and observable, the transmission
zeros are the values of $s$ such that $rank[R(s)]<min(n+m,n+p)$.\end{defn}
\begin{rem}
It is easy to show that for the system considered here, the set of
\{Invariant zeros\} = the set of \{transmission zeros + input decoupling
zeros + output decoupling zeros - input and output decoupling zeros\}. 
\end{rem}
It has been proved in the reference\cite{RefWorks:92} that Definition
\ref{Def:Tzero1} and Definition \ref{Def:Tzero2} are equivalent.
The idea is that the $G(s)$ of a non-degenerate system has a normal
rank of $min(n+m,n+p)$, and therefore the invariant zeros are the
$s$ such that $rank[R(s)]<min(n+m,n+p)$. The condition of controllability
and observability ensures the invariant zeros excludes the input-decoupling
zeros and output-decoupling zeros that are not shown in the $G(s)$.
Then the invariant zeros are indeed the transmission zeros. Then $G(s)$
loses rank at the exactly $s$ at which $R(s)$ loses rank. 
\begin{defn}
The system is minimum phase if all its transmission zeros are in strictly
left-hand plane.
\end{defn}

\section{Squaring-Up Method}

\subsection{Problem Definition}

For the particular interest of adaptive control\cite{RefWorks:90}\cite{RefWorks:36},
we are dealing with a system $\Sigma_{p}=\{A,B,C\}$ that satisfies
following assumptions. \begin{assm}\label{enu:01} The system is
fat and satisfies $n>>m>p$; \end{assm} \begin{assm}The given system
$\{A,B,C\}$ is strictly proper, i.e. $D=0$; \end{assm} \begin{assm}\label{enu:ABC}
$(A,B)$ is a controllable pair, and $(A,C)$ is an observable pair;
\end{assm} \begin{assm}\label{enu:B} $B$ has full column rank,
i.e. $rank(B)=m$; \end{assm} \begin{assm}\label{enu:CB} $rank(CB)=p$
\end{assm}

The observability of $(A,C)$ is not necessarily required for the
following procedure to work\cite{RefWorks:70}. It only serves for
the purpose of analysis simplification. Otherwise, the system will
have output-decoupling zeros that makes $R(s)$ rank deficient according
to Definition \ref{Def:Output Decoup}. We excludes this case since
unobservable modes are not our interest in the context of adaptive
control design. With the observability condition, the only $s$ that
makes $R(s)$ rank deficient is its transmission zeros. 

The goal is to find an augmentation $C_{a}\in\mathbb{R}^{(m-p)\times n}$
such that the system $\{A,B,\bar{C}\}$ is square and minimum phase,
where $\bar{C}^{T}=[C^{T},C_{a}^{T}]$. Assumption \ref{enu:ABC},
\ref{enu:B} and \ref{enu:CB} guarantee $\Sigma_{p}$ is non-degenerate.
Assumption \ref{enu:ABC} guarantee the $\{A,B,C\}$ is the minimal
realization of the system. By Definition \ref{Def:Tzero2}, the number
of transmission zeros should be zero or a finite number. Based on
Definition \ref{Def:Tzero2}, Misra proposes a method to find $C_{a}$
using the technique of pole placement in a special coordinate\cite{RefWorks:70}.
The following section briefly summarizes the steps. For expediency,
we will not distinguish between the term ``rank'' and the term ``normal
rank''.

\subsection{The Existing Method}

The Rosenbrock matrix $R(s)$ can be transformed into a special coordinate
where the controllable states and the uncontrollable states are separated.
\begin{equation}
\tilde{R}(s)=TR(s)T^{T}=\left[\begin{array}{cc|c}
sI_{m}-A_{11} & -A_{12} & -B_{1}\\
-A_{21} & sI_{n-m}-A_{22} & 0\\
\hline C_{11} & C_{12} & 0
\end{array}\right]\label{eq:Tran Rs}
\end{equation}
Since $T$ is an invertible matrix, $rank[\tilde{R}(s)]=rank[R(s)]$
for all $s\in\mathbb{C}$. Assuming a $C_{a}=[C_{21},C_{22}]$ is
found and the augmented system is following: 
\begin{equation}
\tilde{R}_{a}(s)=\left[\begin{array}{cc|c}
sI_{m}-A_{11} & -A_{12} & -B_{1}\\
-A_{21} & sI_{n-m}-A_{22} & 0\\
\hline C_{11} & C_{12} & 0\\
C_{21} & C_{22} & 0
\end{array}\right]\label{eq:Tran Ras}
\end{equation}
Group the column of $\bar{C}$ and denote $C_{1}=\left[\begin{array}{c}
C_{11}\\
C_{21}
\end{array}\right]$ and $C_{2}=\left[\begin{array}{c}
C_{12}\\
C_{22}
\end{array}\right]$. Choose $C_{21}$ such that $C_{1}$ is an invertible matrix. Without
loss of generality, we choose: 
\begin{equation}
C_{21}=null(C_{11}^{T}))\label{eq:C21}
\end{equation}
where $null$ stands for the null space. $C_{21}$ can be made unity.
The following equality holds: {\tiny 
\begin{eqnarray}
 &  & rank[\tilde{R}_{a}(s)]=rank\left(\tilde{R}(s)\left[\begin{array}{cc|c}
I_{m} & -C_{1}^{-1}C_{2} & 0\\
0 & I_{n-m} & 0\\
\hline 0 & 0 & I_{m}
\end{array}\right]\right)\nonumber \\
 &  & =rank\left[\begin{array}{cc|c}
sI_{m}-A_{11} & -sC_{1}^{-1}C_{2}-A_{12}+A_{11}C_{1}^{-1}C_{2} & -B_{1}\\
-A_{21} & sI_{n-m}-A_{22}+A_{21}C_{1}^{-1}C_{2} & 0\\
\hline C_{1} & 0 & 0
\end{array}\right]\label{eq:rank tran}
\end{eqnarray}
}It is easy to show that: 
\begin{eqnarray}
 & rank[\tilde{R}_{a}(s)] & =rank(C_{1})+rank(B_{1})\nonumber \\
 &  & +rank(sI_{n-m}-A_{22}+A_{21}C_{1}^{-1}C_{2})\label{eq:rank cond}
\end{eqnarray}
Since $rank(C_{1})=rank(B_{1})=m$ by design and assumption \ref{enu:B}
and \ref{enu:CB}, $\tilde{R}(s)$ loses rank only if $(sI_{n-m}-A_{22}+A_{21}C_{1}^{-1}C_{2})$
loses rank. From Definition \ref{Def:Tzero2}, the transmission zeros
of the system is exactly the poles of $(A_{22}-A_{21}C_{1}^{-1}C_{2})$.
From Eq.(\ref{eq:Tran Rs}), it is easy to see that the assumption
$(A,B)$ is controllable implies $(A_{22},A_{21})$ is controllable,
which in turn implies $(A_{22},A_{21}C_{1}^{-1})$ is controllable
(since $C_{1}$ is invertible). This implies state feedback technique
can be used on the pair $(A_{22},A_{21}C_{1}^{-1})$ to place the
zeros of the system. The remaining problem is to deal with the fact
that $C_{2}$ is not totally free (since $C_{12}$ is given). Perform
partition on $C_{2}$: 
\begin{equation}
\tilde{C}_{2}=\left[\begin{array}{c}
C_{12}\\
O_{(m-p)\times(n-m)}
\end{array}\right]\; and\;\hat{C}_{2}=\left[\begin{array}{c}
O_{p\times(n-m)}\\
C_{22}
\end{array}\right]
\end{equation}
Correspondingly: 
\begin{equation}
A_{22}-A_{21}C_{1}^{-1}C_{2}=\underset{\tilde{A}_{22}}{\underbrace{A_{22}-A_{21}C_{1}^{-1}\tilde{C}_{2}}}-A_{21}C_{1}^{-1}\hat{C}_{2}
\end{equation}
We only have freedom in designing $\hat{C}_{2}$. That means only
last $(m-p)$ pseudo inputs of $A_{21}C_{1}^{-1}$ are available for
pole placement. Denote $B_{ps}=A_{21}C_{1}^{-1}$ and perform corresponding
partition: 
\begin{equation}
B_{ps}\triangleq A_{21}C_{1}^{-1}\triangleq[B_{ps1},B_{ps2}]
\end{equation}
Now the problem becomes pole placement using feedback on the pair
$(\tilde{A}_{22},B_{ps2})$ where $\tilde{A}_{22}\triangleq A_{22}-A_{21}C_{1}^{-1}\tilde{C}_{2}$.
However, from all above derivation, there is no guarantee that $(\tilde{A}_{22},B_{ps2})$
is controllable. It is found that for some special systems satisfying
all assumptions list above, $(\tilde{A}_{22},B_{ps2})$ can be uncontrollable.
For such system, the existing square-up procedure won't work. Following
context will elaborate the properties of such system and propose the
countermeasure.

\subsection{Presence of Transmission Zeros}

Suppose $(\tilde{A}_{22},B_{ps2})$ is not controllable. Then there
exists a scalar $s_{0}$ and a vector $w_{0}$ such that: 
\begin{equation}
w_{0}^{T}[s_{0}I-\tilde{A}_{22},B_{ps2}]=0
\end{equation}
It follows: 
\begin{eqnarray}
 & w_{0}^{T}s_{0}-w_{0}^{T}\tilde{A}_{22}=0\label{eq:w01}\\
 & w_{0}^{T}B_{ps2}=0\label{eq:w02}
\end{eqnarray}
Substituting the definition of $\tilde{A}_{22}$ transforms Eq.(\ref{eq:w01})
into: 
\begin{equation}
w_{0}^{T}s_{0}-w_{0}^{T}A_{22}+w_{0}^{T}A_{21}C_{1}^{-1}\tilde{C}_{2}=0\label{eq:w01t}
\end{equation}
Also it is noted: 
\begin{equation}
w_{0}^{T}A_{21}C_{1}^{-1}=w_{0}^{T}B_{ps}=w_{0}^{T}[B_{ps1},B_{ps2}]=w_{0}^{T}[B_{ps1},0]\label{eq:Bps1}
\end{equation}
The last equality is true because of Eq.(\ref{eq:w02}). Now let's
examine the form of $C_{1}^{-1}$. With loss of generality, $C_{1}^{-1}$
can be written as: 
\begin{equation}
C_{1}^{-1}=[C_{11}^{\dagger},C_{21}^{T}]\label{eq:C1inv}
\end{equation}
where $C_{11}^{\dagger}$ stands for the right inverse of $C_{11}$.
One can easily verify $C_{1}C_{1}^{-1}=I_{m}$ using the facts: 
\begin{equation}
\begin{cases}
 & C_{11}C_{21}^{T}=O_{p\times(m-p)}\\
 & C_{11}C_{11}^{\dagger}=I_{p}\\
 & C_{21}C_{11}^{\dagger}=O_{(m-p)\times p}\\
 & C_{21}C_{21}^{T}=I_{m-p}
\end{cases}
\end{equation}
Using Eq.(\ref{eq:C1inv}), $B_{ps1}$ and $B_{ps2}$ can be rewritten
as: 
\begin{equation}
B_{ps1}=A_{21}C_{11}^{\dagger}\qquad and\qquad B_{ps2}=A_{21}C_{21}^{T}\label{eq:Bps}
\end{equation}
Eq.(\ref{eq:Bps1}) and Eq.(\ref{eq:Bps}) can transform Eq.(\ref{eq:w01t})
into: 
\begin{equation}
w_{0}^{T}s_{0}-w_{0}^{T}A_{22}+w_{0}^{T}A_{21}C_{11}^{\dagger}C_{12}=0\label{eq:rankdef}
\end{equation}
Now we will examine what Eq.(\ref{eq:rankdef}) implies of the original
system (\ref{eq:Tran Rs}). Following equality takes place: {\tiny 
\begin{eqnarray}
 &  & rank[\tilde{R}(s)]=rank\left(\tilde{R}(s)\left[\begin{array}{cc|c}
I_{m} & -C_{11}^{\dagger}C_{12} & 0\\
0 & I_{n-m} & 0\\
\hline 0 & 0 & I_{m}
\end{array}\right]\right)\nonumber \\
 &  & =rank\left[\begin{array}{cc|c}
sI_{m}-A_{11} & -sC_{11}^{\dagger}C_{12}-A_{12}+A_{11}C_{11}^{\dagger}C_{12} & -B_{1}\\
-A_{21} & sI_{n-m}-A_{22}+A_{21}C_{11}^{\dagger}C_{12} & 0\\
\hline C_{11} & 0 & 0
\end{array}\right]
\end{eqnarray}
}Similar to Eq.(\ref{eq:rank cond}), the rank of $\tilde{R}(s)$
fully depends on $C_{11}$, $B_{1}$ and $sI_{n-m}-A_{22}+A_{21}C_{11}^{\dagger}C_{12}$:
\begin{eqnarray}
 & rank[\tilde{R}(s)] & =rank(C_{11})+rank(B_{1})\nonumber \\
 &  & +rank(sI_{n-m}-A_{22}+A_{21}C_{11}^{\dagger}C_{12})\label{eq:Ori rank cond}
\end{eqnarray}
Eq.(\ref{eq:rankdef}) says there exists a $s_{0}$ such that $sI_{n-m}-A_{22}+A_{21}C_{11}^{\dagger}C_{12}$
loses rank. Eq.(\ref{eq:Ori rank cond}) says such $s_{0}$ will make
the original system $\tilde{R}(s)$ loses rank. By definition \ref{Def:Tzero2},
the system $R(s)$ has a transmission zero at $s_{0}$. It can be
concluded now that given assumption \ref{enu:01} to \ref{enu:CB},
the only case Misra's method can not solve is the case when the given
system already posses a transmission zero. Comparing Eq.(\ref{eq:Tran Rs})
and Eq.(\ref{eq:Tran Ras}), it is easy to see that any $s_{0}$ that
makes $\tilde{R}(s)$ lose rank will also make $\tilde{R_{a}}(s)$
loses rank. In other words, any transmission zeros of the given system
will become the transmission zeros of the squared-up system. That
is one important limitation of the square-up procedure.

\subsection{Countermeasure\label{sub:Countermeasure}}

The countermeasure is following. The above derivation can be reversed
and the sufficient condition argument is true, i.e. the transmission
zeros of the given system is indeed the uncontrollable mode of the
pair $(\tilde{A}_{22},B_{ps2})$. Even if uncontrollable modes exist,
other controllable modes can be placed in the strictly left-hand plane
using the remaining feedback action. As a result, it can be concluded
that if the given system satisfies one additional condition: \begin{assm}
\label{enu:Tran}The system has only nonminimum-phase transmission
zeros. \end{assm} the pair $(\tilde{A}_{22},B_{ps2})$ is stabilizable
and we are still be able to design a $C_{a}$ such that the squared-up
system is minimum phase. LQR technique is immediately available for
such problem. We summarized our improved method as following: \begin{stepwise}
Check if the given system satisfies all assumption \ref{enu:01} to
\ref{enu:Tran};\end{stepwise} \begin{stepwise} Transform it into
a controllable canonical form as in Eq.(\ref{eq:Tran Rs}); \end{stepwise}
\begin{stepwise} Find $C_{21}$ using Eq.(\ref{eq:C21}); \end{stepwise}
\begin{stepwise} Calculate the stabilizable pair $(\tilde{A}_{22},B_{ps2})$;
\end{stepwise} \begin{stepwise} Perform LQR technique on $(\tilde{A}_{22},B_{ps2})$
to find $C_{22}$; \end{stepwise} \begin{stepwise} Augment $C$
and transform the system back to its original coordinate. \end{stepwise}

\section{Example\label{sec:Example}}

Following context gives an example of a MIMO system with a transmission
zero and the results of our squaring-up procedure. It is a linearized
model for the lateral dynamics of Boeing 747-100 aircraft. We transposed
the system for the illustration of a fat system. {\tiny 
\begin{equation}
R(s)=\left[\begin{array}{cccc|ccc}
s+0.0605 & 0.0015 & -0.0011 & 0 & 0 & 0 & 0\\
0 & s+0.4603 & 0.0208 & 1 & -1 & 0 & 0\\
871 & -0.28 & s+0.141 & 0 & 0 & -1 & 0\\
32.3 & 0 & 0 & 0 & 0 & 0 & -1\\
\hline 0 & -0.1860 & 0.0061 & 0 & 0 & 0 & 0\\
4.0380 & 0.1 & -0.4419 & 0 & 0 & 0 & 0
\end{array}\right]
\end{equation}
}A quick check can confirm the given system satisfies assumption \ref{enu:01}
to \ref{enu:CB}, and has an transmission zero at $-0.0511$. Coordinate
transformation $T=[B^{T},(null(B^{T}))^{T}]$gives: {\tiny 
\begin{equation}
\tilde{R}(s)=\left[\begin{array}{cccc|ccc}
s+0.4603 & 0.0208 & -1 & 0 & -1 & 0 & 0\\
-0.28 & s+0.1410 & 0 & -871 & 0 & -1 & 0\\
0 & 0 & 0 & -32.3 & 0 & 0 & -1\\
-0.0015 & 0.0011 & 0 & 0.0605 & 0 & 0 & 0\\
\hline -0.1860 & 0.0061 & 0 & 0 & 0 & 0 & 0\\
0.1 & -0.4419 & 0 & -4.038 & 0 & 0 & 0
\end{array}\right]
\end{equation}
}Quick examination sees: 
\begin{equation}
A_{22}=-0.0605\quad and\quad A_{21}=[0.0015,-0.0011,0]
\end{equation}
By our choice of $C_{21}$, $C_{1}$ and $C_{1}^{-1}$ is given: 
\begin{eqnarray}
 & C_{1}=\left[\begin{array}{ccc}
-0.186 & 0.0061 & 0\\
0.1 & -0.4419 & 0\\
0 & 0 & 1
\end{array}\right]\\
 & C_{1}^{-1}=\left[\begin{array}{ccc}
-5.4171 & -0.0744 & 0\\
-1.2263 & -2.2796 & 0\\
0 & 0 & 1
\end{array}\right]
\end{eqnarray}
And $C_{21}$ and $\tilde{C}_{2}$ are: 
\begin{equation}
C_{21}=[0,0,1]\quad and\quad\tilde{C}_{2}=\left[\begin{array}{c}
0\\
-4.038\\
0
\end{array}\right]
\end{equation}
Then $\tilde{A}_{22}=A_{22}-A_{21}C_{1}^{-1}\tilde{C}_{2}$ and $B_{ps2}=A_{21}C_{21}^{T}$
gives: 
\begin{equation}
\tilde{A}_{22}=-0.0511\quad and\quad B_{ps2}=0
\end{equation}
It is verified that the pair $(\tilde{A}_{22},B_{ps2})$ is uncontrollable
and the uncontrollable mode is exactly the transmission zero of the
given system. Simply put $\hat{C}_{2}=O_{m\times(n-m)}$ and transform
the $C_{a}$ back into the original coordinate. The augmented system
will become: {\tiny 
\begin{equation}
R_{a}(s)=\left[\begin{array}{cccc|ccc}
s+0.0605 & 0.0015 & -0.0011 & 0 & 0 & 0 & 0\\
0 & s+0.4603 & 0.0208 & 1 & -1 & 0 & 0\\
871 & -0.28 & s+0.141 & 0 & 0 & -1 & 0\\
32.3 & 0 & 0 & 0 & 0 & 0 & -1\\
\hline 0 & -0.1860 & 0.0061 & 0 & 0 & 0 & 0\\
4.0380 & 0.1 & -0.4419 & 0 & 0 & 0 & 0\\
0 & 0 & 0 & 1 & 0 & 0 & 0
\end{array}\right]
\end{equation}
}The last row is the designed pseudo-output. Quick examination will
verify the augmented system $R_{a}(s)$ has only one transmission
zero at $-0.0511$.

\section{Conclusions}

This paper proves that the uncontrollable modes in the existing squaring-up
method are actually the transmission zeros of the given system. In
other words, we are not able to move the locations of the existing
transmission zeros using the proposed method. Systems with minimum
phase zeros has a stabilizable pair in the transformed coordinate
and therefore can be squared-up using LQR technique. It is noted that
by transposing the system, the proposed method can be applied on a
tall system with more outputs than inputs.

\bibliographystyle{IEEEtran}
\bibliography{Ref/VFA_Refworks_v4}

\end{document}